\documentclass{amsart}

\usepackage{amssymb}

\title{FUNCTIONAL REPESENTATION OF SUBSTITUTION ALGEBRAS }

\author{Norman Feldman}

\address{Mathematics Department\\Sonoma State University\\Rohnert Park, CA 94928, USA}

\email{norm.feldman@sonoma.edu}

\newtheorem{Theorem}{Theorem}[section]

\theoremstyle{definition}
\newtheorem{Definition}[Theorem]{Definition}

\begin{document}

\begin{abstract}
We show that the class of representable substitution algebras is  characterized by a set of universal first order sentences. In addition, it is shown that a necessary and sufficient condition for a substitution algebra to be representable is that it is embeddable in a substitution algebra in which elements are distinguished. Furthermore, conditions in terms of neat embeddings are shown to be equivalent to representability.
\end{abstract}

\maketitle

\section{Introduction}
In [1], the problem of axiomatizing polonomial substitution algebras was considered. It was shown that a set of first order sentences and a non-first order condition of local finitness characterizes the class of isomorphs of these algebras. From Theorem 3.1 in [1], it follows that every locally finite substitution algebra is isomorphic to a function substitution algebra and hence, is representable.

In this section we recall some of the definitions from [1].

Let $\alpha$ be an ordinal and let $F_\alpha(U)$ be the set of functions from $U^{\alpha}$, the set of $\alpha$-sequences of members of $U$, to $U$, a non-empty set. 

\begin{Definition}
For $s \in U^{\alpha}$, let $s\langle \kappa, x\rangle \in U^{\alpha}$ be defined by $s\langle \kappa, x\rangle_\lambda = s_\lambda$ if $\lambda \neq \kappa$ and  $s\langle \kappa, x\rangle_\lambda = x$ if $ \lambda = \kappa$.
\end{Definition}

\begin{Definition}
For $\kappa<\alpha$, define a binary operation on  $F_\alpha(U)$ by 
\[
(f*_kg)(s)=g(s\langle k,f(s)\rangle)
\]
 for $f, g\in F_\alpha(U)$.
\end{Definition}

Let $V_\kappa\in F_\alpha(U)$ be defined by $V_\kappa\left(s\right)=s_\kappa$ for $s\in U^{\alpha}$.  

\begin{Definition}
An algebra of the form $\mathfrak{A}=\langle A, *_\kappa, V_\kappa\rangle_{\kappa<\alpha}$ where $A\subseteq F_\alpha(U)$, $A$ is closed under $*_\kappa$, and $V_\kappa\in A$ for $\kappa<\alpha$, is a function substitution algebra of dimension $\alpha$ ($FSA_\alpha$). If $A=F_\alpha(U)$, then $\mathfrak{A}$ is said to be the full function substitution algebra of dimension $\alpha$ with base $U$.
\end{Definition}

\begin{Definition}
An algebra $\mathfrak{A}=\langle A, *_\kappa, v_\kappa\rangle_{\kappa<\alpha}$ satisfying the following first order axiom schema is called a substitution algebra of dimension $\alpha$  ($SA_\alpha$):

(1) For all $x$, $x*_\kappa v_\kappa=x$.

(2) For all $x$, $x*_\kappa v_\lambda=v_\lambda$ for $\kappa \neq \lambda$.

(3) For all $x$, $v_\kappa*_\kappa x=x$.

(4)  For all $x$ and $z$, if $w*_\lambda x=x$, for all $w$, then $x*_\lambda \left(v_\lambda *_\kappa z \right)=x*_\lambda \left(x *_\kappa z \right)$.

(5) For all $x$, $y$, and $z$, $\left(x *_\kappa y \right)*_\kappa z=x *_\kappa \left(y*_\kappa z \right)$.

(6) For all $x$, $y$, and $z$, if $w*_\lambda x=x$, for all $w$, then 
\[
x *_\kappa \left(y*_\lambda z \right) = \left(x *_\kappa y \right)*_\lambda \left(x*_\kappa z \right) \textrm{, for } \kappa \neq \lambda.
\]
\end{Definition}

It is easily shown that every $FSA_\alpha$ is an $SA_\alpha$.

In [1] it was shown that the class of $SA_\alpha$'s can be axiomatized by a set of universal sentences and hence, a subalgebra of a $SA_\alpha$ is an $SA_\alpha$.

\begin{Definition}
Let $\mathfrak{A}=\langle A, *_\kappa, v_\kappa\rangle_{\kappa<\alpha}$ be an $SA_\alpha$. For $x \in A$, let
\[
\Delta _\mathfrak{A}x=\{\kappa:a*_\kappa x\neq x, \textrm{ for some } a\in A \}.
\]
$\Delta _\mathfrak{A}x$ is called the dimension set of $x$.
\end{Definition}

We usually drop the subscript and write $\Delta x$ for $\Delta _\mathfrak{A}x$ when there is no ambiguity.

\begin{Definition}
A substitution algebra $\mathfrak{A}=\langle A, *_\kappa, v_\kappa\rangle_{\kappa<\alpha}$ is locally finite
if $\Delta x$ is finite for all $x \in A$.
\end{Definition}

For an arbitrary algebra $\mathfrak{B}$, the polynomial substitution algebra over $\mathfrak{B}$ of dimension $\alpha$ is a locally finite $SA_\alpha$. It is, in fact, an $FSA_\alpha$. In [1]. it was shown that every locally finite $SA_\alpha$ is isomorphic to some polynomial substitution algebra. The question arises: Is every $SA_\alpha$ isomorphic to an $FSA_\alpha$? Call such an algebra a representable substitution algebra ($RSA_\alpha$).

In section 3, we show that the class of $RSA_\alpha$ is at least characterized by a set of first order universal sentences; that is, the class  $RSA_\alpha$ is a $UC_\Delta$.

\begin{Definition}
Elements are distinguished in an $SA_\alpha$ if for all $\kappa < \alpha$, $x\neq y$ implies $c*_\kappa x \neq c*_\kappa y$ for some $c$ with $\Delta c=0$.
\end{Definition}

In section 4, we show that an $SA_\alpha$ is an $RSA_\alpha$ iff it is embeddable in an $SA_\alpha$ in which elements are distinguished. At this time, we do not know if there is an $SA_\alpha$ which is not an $RSA_\alpha$.

Unless otherwise stated, the universe of a substitution algebra will be denoted by the Roman letter corresponding to the Gothic letter used to denote the algebra. For example, the universes of $\mathfrak{A}$, $\mathfrak{B}_i$,and $\mathfrak{C}'$ are $A$, $B_i$, and $C'$ respectively.

Filters will be used in several places. Let $\mathcal{F}$ be a filter on a set $I$ and let $X_i$ be a set for $x\in I$. Let $\thicksim$ be the equivalence relation on the product $W=\prod \langle X_i:i\in I \rangle$ defined by $x \thicksim y$ iff $\{i:i \in I \textrm{ and } x_i = y_i\} \in \mathcal{F}$ for $x, y \in W$.

Throughout the paper the equivalence class containing $x$ is denoted by $[x]$.

\section{$\Gamma$-homomorphisms}

\begin{Definition}
Let $\mathfrak{A}$ be an $SA_\alpha$. For a subset $\Gamma$ of $\alpha$, let 
\[
Z(\Gamma)=\{x:x \in A \textrm{ and } \Delta x \cap \Gamma=0\}.
\]
\end{Definition}

It follows that $Z(\alpha)=\{x:x \in A \textrm{ and } \Delta x=0\}$. In this case, let $Z=Z(\alpha)$.

In [1], the notion of generalized substitutions was defined by
\[
s*_{(\Gamma)}a=s_{\kappa_{n-1}} *_{\kappa_{n-1}}(\cdots(s_{\kappa_1} *_{\kappa_1}\left(s_{\kappa_0}*_{\kappa_0}a))\cdots\right)
\]
where $s\in Z^\alpha$, $a\in A$ and $\Gamma=\{\kappa_0,\cdots,\kappa_{n-1}\}$ with $\kappa_0<\cdots<\kappa_{n-1}$.

We can let $s\in A^\alpha$ in this definition. With suitable restrictions on $s$, we obtain theorems analogous to Theorems 4.1  and 4.2 of [1]. We state below the parts of these theorems needed for this paper. The proofs are almost identical to their counterparts in [1]. The following theorem shows that a generalized substitution, $s*_{(\Gamma)}a$ is independent of the order of the elements in $\Gamma$.

\begin{Theorem}
If $s\in Z(\Gamma)^\alpha$, $a\in A$, and $\Gamma=\{\lambda_0,\cdots ,\lambda_{n-1}\}$, then 
\[
s*_{(\Gamma)}a=s_{\lambda_{n-1}} *_{\lambda_{n-1}}(\cdots(s_{\lambda_1} *_{\lambda_1}(s_{\lambda_0}*_{\lambda_0}a))\cdots).
\]
\end{Theorem}

\begin{Theorem}
Let $\Gamma$ and $\Sigma$ be finite subsets of $\alpha$.

(i) If $s\in Z(\Gamma)^\alpha$ and $\kappa \in  \Gamma $, then $s*_{(\Gamma)}a=s*_{(\Gamma-\{\kappa\})}(s_\kappa*_\kappa a)$.

(ii) If $s\in Z(\Gamma)^\alpha$, then $\Delta (s*_{(\Gamma)} a)\subseteq \bigcup\{\Delta s_\kappa: \kappa \in\Gamma\}\cup(\Delta a-\Gamma)$.

(iii) If $s\in Z(\Gamma)^\alpha$, $x \in Z(\Gamma)$, and $\kappa \in \Gamma$, then $s\langle \kappa,x\rangle*_{(\Gamma)}a=x*_\kappa(s*_{(\Gamma-\{\kappa\})}a)$.

(iv) If $s\in Z(\Gamma\cup\Sigma\cup\{\kappa\})^\alpha$, and $\kappa \notin \Sigma$, then
\[
(s*_{(\Gamma)}a)*_\kappa(s*_{(\Sigma)}b)=s*_{(\Gamma \cap \Sigma)}((s*_{(\Gamma-\Sigma)}a)*_\kappa (s*_{(\Sigma-\Gamma)}b)).
\]
\end{Theorem}

\begin{Definition}
Let $\mathfrak{A}=\langle A,*_\kappa,v_\kappa \rangle_{\kappa<\alpha}$ and $\mathfrak{A}'=\langle A',*'_\kappa,v'_\kappa \rangle_{\kappa<\alpha}$ be $SA_\alpha$'s. and let $\Gamma\subseteq \alpha$. $\phi$ is a $\Gamma$-homomorphism from $\mathfrak{A}$ into $\mathfrak{A}'$ if $\phi:A\to A'$ and for all $\kappa \in \Gamma$ and all $a,b \in A$, $\phi(a*_\kappa b)=\phi(a)*'_\kappa \phi(b)$ and $\phi(v_\kappa)=v_\kappa'$.
\end{Definition}

\begin{Theorem}
Let $\mathfrak{A}$ be an $SA_\alpha$ and let $\Gamma$ be a finite subset of $\alpha$. Define a function $\phi_\Gamma:A\to F_\alpha(Z(\Gamma))$ by $\phi_\Gamma(a)(s)=s*_{(\Gamma)}a$. Then $\phi_\Gamma$ is a $\Gamma$-homomorphism from $\mathfrak{A}$ to the full $FSA_\alpha$ with base $Z(\Gamma)$.

\begin{proof}
First note that Theorem 2.3 (ii) implies that $s*_{(\Gamma)} a \in Z(\Gamma)$ for $s \in Z(\Gamma)^\alpha$ and $s \in A$. 

It is easily checked that $\phi_\Gamma(v_\kappa)=V_\kappa$ for $\kappa \in \Gamma$ where $V_\kappa(s)=s_\kappa$ for $s\in Z(\Gamma)^\alpha$.

Let $\kappa \in \Gamma$, $a, b \in A$, and $s  \in Z(\Gamma)^\alpha$. We then have

 \begin{eqnarray*}
(\phi_\Gamma(a)*_\kappa \phi_\Gamma(b))(s) & =& \phi_\Gamma(b)(s\langle \kappa, \phi_\Gamma(a)(s)\rangle)\\
 & = & s\langle \kappa,s*_{(\Gamma)}a\rangle*_{(\Gamma)}b\\
 & = & (s*_{(\Gamma)}a)*_\kappa(s*_{(\Gamma-\{\kappa\})}b)      \textrm{   by Theorem 2.3 (ii), (iii)}\\
 & = & s*_{(\Gamma-\{\kappa\})}((s*_\kappa a)*_\kappa b)  \textrm{   by Theorem 2.3 (iv)}\\
 & = & s*_{(\Gamma-\{\kappa\})}(s_\kappa *_\kappa(a*_\kappa b))\\
 & = & s_{(\Gamma)}(a*_\kappa b)     \textrm{   by Theorem 2.3 (i)}\\
 & = & \phi_{(\Gamma)}(a*_\kappa b)(s).
\end{eqnarray*}
\end{proof}
\end{Theorem}

\section{A representation theorem and ultraproducts}

\begin{Definition}
Let $\mathfrak{A}$ be an $SA_\alpha$. Elements are strongly distinguished in $\mathfrak{A}$ if, for all $s\in Z^\alpha$, there is a finite subset $\Sigma$ of $\alpha$, such that $s*_{(\Sigma)}a=s*_{(\Sigma)}b$, then $a=b$.
\end{Definition}

\begin{Theorem}
Every $SA_\alpha$ in which elements are strongly distinguished is an $RSA_\alpha$.

\begin{proof}

Let $\mathfrak{A}$ be an $SA_\alpha$ in which elements are strongly distinguished. Let $J$ be the collection of all finite subsets of $\alpha$ and let
\[
P_\Sigma=\{\Gamma:\Gamma \in J \textrm{ and } \Gamma \supseteq\Sigma\}
\]
for $\Sigma \in J$. 
Since $P_\Sigma \cap P_\Gamma=P_{\Sigma \cup \Gamma}$, it follows that
\[
\mathcal{F}=\{Q:Q \subseteq J \textrm{ and for some } \Sigma \in J, Q\supseteq P_\Sigma\}
\]
is a filter on $J$.

Let $W=\prod\langle Z(\Gamma):\Gamma \in J\rangle$, and let $U$ be the collection of all equivalence classes $[f]$ for $f \in W$. $U$ is the base for the $FSA_\alpha$ we seek.

For $X \in U$, let
\[
K(X)=\{f:f \in X \textrm{ and for some } x\in Z, f_\Gamma=x \textrm{ for all } \Gamma \in J\}.
\]
It is easily seen that $K(X)$ has at most one member. Let $c:U\to W$ be a choice function such that for all $X \in U$, $c(X) \in K(X)$ if $K(X) \ne \oslash$ and $c(X) \in X$ if $K(X)=\oslash$.

For each $s \in U^\alpha$ and $\Gamma \in J$, define $s^\Gamma \in Z(\Gamma)^\alpha$ by $(s^\Gamma)_\lambda=c(s_\lambda)$ for $\lambda <\alpha$. For $\Gamma \in J$, let $\phi_\Gamma$ be the  $\Gamma$- homomorphism of Theorem 2.5. Define a function $ \phi: A \to F_\alpha(U)$ by $\phi(a)(s)=[\langle \phi_\Gamma(a)(s^\Gamma):\Gamma \in J \rangle]$ where $a \in A$ and $s \in U^\alpha$. We shall show that $\phi$ is an isomorphism from $\mathfrak{A}$ to the full $FSA_\alpha$ with base $U$. The result follows immediately from this.

First we show that $\phi$ is one-one. Suppose that $\phi(a)=\phi(b)$ where $a,b \in A$.  Therefore, for all $s\in U^\alpha$, $\phi(a)(s)=\phi(b)(s)$ and hence,
\[
\{\Gamma:\Gamma \in J \textrm{ and } s^\Gamma*_{(\Gamma)}a=s^\Gamma*_{(\Gamma)}b\} \in \mathcal{F}.
\]
From the definition of $\mathcal{F}$, it follows that 
\begin{equation}
\textrm{ for all } s \in U^\alpha, \textrm{ there is a } \Sigma \in J \textrm{ such that } s^\Sigma*_{(\Sigma)}a =s^\Sigma*_{(\Sigma)}b.
\end{equation}

Since elements are strongly distinguished in $\mathfrak{A}$, to show that $a=b$,, it suffices to show that 
\begin{equation}
\textrm{for all }t \in Z^\alpha, \textrm{ there is a } \Sigma \in J \textrm{ such that }t*_{(\Sigma)}a=t*_{(\Sigma)}b.
\end{equation}

To this end, let $t \in Z^\alpha$ and define $s\in U^\alpha$  by $s_\lambda=[\langle t_\lambda: \Gamma \in J\rangle]$ for $\lambda<\alpha$. By (1), there is a $\Sigma \in J$ such that $s^\Sigma*_{(\Sigma)}a=s^\Sigma*_{(\Sigma)}b$.  Since $\langle t_\lambda:\Gamma \in J\rangle \in K(s_\lambda)$ we have $c(s_\lambda)=\langle t_\lambda:\Gamma \in J\rangle$. Therefore,, for all $\lambda<\alpha$, $(s^\Sigma)_\lambda=c(s_\lambda)^\Gamma=t_\lambda$. Hence, $s^\Sigma=t $ and $t*_{(\Sigma)}a=t*_{(\Sigma)}b$. We therefore have (2).

Next we show that $\phi$ is a homomorphism from $\mathfrak{A}$ to to the full $FSA_\alpha$ with base $U$.  We need the following easily verified fact: 
\begin{equation}
\textrm{For all } s \in U^\alpha \textrm{and } X\in U, s\langle \kappa,x\rangle^\Gamma=s^\Gamma\langle \kappa,c(X)_\Gamma\rangle \textrm{.}
\end{equation}

Let $a,b \in A$ and $s \in U^\alpha$. Let
\[
M=\{\Gamma:\Gamma \in J \textrm{ and } \phi_\Gamma(a*_\kappa b)(s^\Gamma)=(\phi_\Gamma(a)*_\kappa\phi_\Gamma(b))(s^\Gamma)\}.
\]

Since $\phi_\Gamma$ is a $\Gamma$-homomorphism, $M\supseteq P_{\{\kappa\}}$ and hence, $M \in \mathcal{F}$. We then have

\begin{eqnarray*}
\phi(a*_\kappa b)(s) & = & [\langle \phi_\Gamma(a*_\kappa b)(s^\Gamma):\Gamma \in J \rangle] \\
& = & [\langle (\phi_\Gamma(a)*_\kappa \phi_\Gamma(b))(s^\Gamma):\Gamma \in J \rangle]  \textrm{    since M } \in \mathcal{F} \\
& = & [\langle \phi_\Gamma(b)(s^\Gamma\langle \kappa,\phi_\Gamma(a)s^\Gamma)\rangle):\Gamma \in J\rangle].
\end{eqnarray*}

On the other hand, 
\begin{eqnarray*}
(\phi(a)*_\kappa\phi(b))(s) & = & \phi(b)(s\langle \kappa,\phi(a)(s)\rangle)\\
& = & [\langle\phi_\Gamma(b)(s\langle \kappa,\phi(a)(s)\rangle^\Gamma):\Gamma \in J\rangle]\\
& = & [\langle \phi_\Gamma(b)(s^\Gamma\langle \kappa,c(\phi(a)(s)_\Gamma\rangle):\Gamma \in J\rangle].
\end{eqnarray*}

Let
\[
R=\{\Gamma:\Gamma \in J \textrm{ and } \phi_\Gamma(a)(s^\Gamma)=c(\phi(a)(s))_\Gamma\}.
\]
It is easily seen that $R \in \mathcal{F}$. It follows that $\phi(a*_\kappa b)(s)=(\phi(a)*_\kappa \phi(b))(s)$.

By a similar type of argument, it can be shown that $\phi(v_\kappa)=V_\kappa$ where $V_\kappa(s)=s_\kappa$ for all $s \in U^\alpha$.
\end{proof}
\end{Theorem}

\begin{Theorem}
Elements are strongly distinguished in every full $FSA_\alpha$.

\begin{proof}
Let $\mathfrak{A}$ be a full $FSA_\alpha$ with base $U$. Suppose that $f, g \in F_\alpha(U)$ such that $f \ne g$. Therefore, $f(s) \ne g(s)$ for some $s \in U^\alpha$. Define $t \in Z^\alpha$ by $t_\lambda(r) = s_\lambda$ for $r \in U^\alpha$. A direct computation shows that for all finite subsets $\Sigma$ of $\alpha$, $(t*_{(\Sigma)}f)(s) = f(s)$ and $(t*_{(\Sigma)}g)(s) = g(s)$ and hence, $t*_{(\Sigma)}g \ne t*_{(\Sigma)}f$. Therefore, we have shown that if $f \ne g$, then there is a $t \in Z^\alpha$ such that $t*_{(\Sigma)}g \ne t*_{(\Sigma)}f$ for all finite subsets $\Sigma$ of $\alpha$.
\end{proof}
\end{Theorem}

\begin{Theorem}
An $SA_\alpha$ is representable iff it is embeddable in an $SA_\alpha$ in which elements are strongly distinguished.
\begin{proof}
If  an $SA_\alpha$ is representable, then it is isomorphic to an $FSA_\alpha$ and, since every $FSA_\alpha$ is a subalgebra of a full $FSA_\alpha$, by Theorem 3.3, it is embeddable in an $SA_\alpha$ in which elements are strongly distinguished.

The other implication follows from Theorem 3.2.
\end{proof}
\end{Theorem}

We use the next theorem to show that an ultraproduct of $FSA_\alpha$'s is representable.

\begin{Theorem}
Let $\mathcal{F}$ be an ultrafilter on a set $I$ and for all $i\in I$, let $\mathfrak{A_i}$ be an $SA_\alpha$ in which elements are strongly distinguished. Then elements are strongly distinguished in $\mathfrak{B}=\prod \langle \mathfrak{A_i}:i\in I \rangle/\mathcal{F}$, the ultra product of the $\mathfrak{A_i}$.

\begin{proof}
The universe $B$ of $\mathfrak{B}$ is the collection of equivalence classes $[a]$ where $a \in \prod\langle A_i:i \in I\rangle$. Let $a,b \in W$ such that $[a] \ne [b]$.  Let $P=\{i:i\in I \textrm{ and } a_i \ne b_i\}$. Since $\mathcal{F}$ is an ultrafilter, $P\in \mathcal{F}$. Let $Z_i = \{x:x\in A_i \textrm{ and } \Delta x=0\}$ and $Z = \{x:x \in B \textrm{ and }  \Delta x = 0\}$. For all $i \in I$, elements are strongly distinguished in $\mathfrak{A_i}$. Therefore, for all $ i \in P$, there is a $t^i \in Z^\alpha_i$ such that $t^i*_{(\Sigma)}a_i \ne t^i*_{(\Sigma)}b_i$ for all finite subsets $ \Sigma$ of $\alpha$. For $i \in I-P$, let $t^i$ be any member of $Z^\alpha_i$. Define $s \in Z^\alpha$ by letting $s_\lambda = [\langle (t^i)_\lambda:i \in I \rangle]$ for $\lambda < \alpha$. 

Let $\Sigma$ be a finite subset of $\alpha$. A direct computation shows that for all $x \in W$, $s*_{(\Sigma)}[x]=[\langle t^i*_{(\Sigma)}x_i:i \in I\rangle]$. Let $Q = \{i:i \in I \textrm{ and } t^i*_{(\Sigma)}a_i \ne t^i*_{(\Sigma)} b_i\}$.

Since $P \in \mathcal{F}$ and $Q\supseteq P$, we have $Q \in \mathcal{F}$. Therefore, 
\begin{eqnarray*}
s*_{(\Sigma)}[a] & = & [\langle t^i *_{(\Sigma)} a_i:i \in I\rangle]\\
& \ne & [\langle t^i *_{(\Sigma)}b_i;i\in I\rangle] \\
& = & s*_{(\Sigma)}[b].
\end{eqnarray*}

\end{proof}
\end{Theorem}

\begin{Theorem}
An ultraproduct of $RSA_\alpha$'s is an $RSA_\alpha$.
\begin{proof}
For $i \in I$ let $\mathfrak{A_i}$ be an $RSA_\alpha$ and let $\mathcal{F}$ be an ultrafilter  on $I$. Each $\mathfrak{A_i}$ is isomorphic to some $FSA_\alpha$ $\mathfrak{A_i}'$ with base $U_i$. $\mathfrak{A_i}'$ is a subalgebra of the full $FSA_\alpha$ $ \mathfrak{B_i}$ with base $U_i$. $\mathfrak{A}=\prod \langle \mathfrak{A_i}:i \in I\rangle/\mathcal{F}$ is isomorphic to a subalgebra of $\mathfrak{B}=\prod \langle \mathfrak{B_i}:i \in I\rangle/\mathcal{F}$. By Theorems 3.3, 3.5, and 3.2, $\mathfrak{B}$ is an $RSA_\alpha$ and therefore, $\mathfrak{A}$ is an $RSA_\alpha$.
\end{proof}
\end{Theorem}

\begin{Theorem}
The class of $RSA_\alpha$'s is a $UC_\Delta$.
\begin{proof}
The class of $RSA_\alpha$'s is closed under isomorphism, subalgebra, and by Theorem 3.6, ultraproducts.  By a theorem of Los [4], the class of $RSA_\alpha$'s is a $UC_\Delta$.
\end{proof}
\end{Theorem}

\section{A condition equivalent to representability}

It is easily seen that if elements are strongly distinguished in an $SA_\alpha$ $\mathfrak{A}$ , then elements are distinguished in $\mathfrak{A}$. Therefore, by Theorem 3.3, elements are distinguished in a full $FSA_\alpha$.

\begin{Theorem}
An $SA_\alpha$ is an $RSA_\alpha$ iff it is emeddable in an $SA_\alpha$ in which elements are distinguished.
\begin{proof}
Every $FSA_\alpha$ is a subalgebra of a full  $FSA_\alpha$ and elements are distinguished in every full   $FSA_\alpha$. Therefore, every representable  $SA_\alpha$ is embeddable in an  $SA_\alpha$ in which elements are distinguished.

For the converse, let $\mathfrak{A}'$ be an $SA_\alpha$ which is embeddable in $\mathfrak{A}$ in which elements are distinguished. To show that $\mathfrak{A}'$ is representable, it suffices to show that $\mathfrak{A}$ is representable.

Let $J$ be the collection of finite subsets of $\alpha$. For $\Gamma \in J$, let $\mathfrak{A}_\Gamma$ be the full $FSA_\alpha$ with base  $Z(\Gamma)$. For all $\Gamma \in J$, $\phi_\Gamma$, the function defined in Theorem 2.5, is a $\Gamma$-homomorphism from $\mathfrak{A}$ into $\mathfrak{A}_\Gamma$. We show that $\phi_\Gamma$ is one-one. Suppose that $\phi_\Gamma(a)=\phi_\Gamma(b)$ where $a, b \in A$. Hence, for all $s \in Z(\Gamma)^\alpha$, $\phi_\Gamma(a)(s)=\phi_\Gamma(b)(s)$; that is, $s*_{(\Gamma)}a=s*_{(\Gamma)}b$. Since $Z \subseteq Z(\Gamma)$, $s*_{(\Gamma)}a=s*_{(\Gamma)}b$ for all $s \in Z^\alpha$. By Theorem 5.5 of [1], $a=b$ since elements are distinguished in $\mathfrak{A}$. 

For $\lambda < \alpha$, let $R_\lambda=\{\Gamma:\Gamma \in J \textrm{ and } \lambda \in \Gamma\}$. The collection of all $R_\lambda$ for $\lambda < \alpha$, has the finite intersection property and hence, there is an ultrafilter $\mathcal{F}$ on $J$ such that $R_\lambda \in \mathcal{F}$ for all $\lambda < \alpha$. Let $\mathfrak{B}=\prod\langle  \mathfrak{A_\Gamma}:\Gamma \in J\rangle/\mathcal{F}$. Define a map $\psi:A \to B$ as follows: $\psi(a)=[\langle \phi_\Gamma(a):\Gamma \in J\rangle]$ for $a \in A$. A straightforward verification shows that $\psi$ is an isomorphism from $\mathfrak{A}$ into  $\mathfrak{B}$. By Theorem 3.6, $\mathfrak{B}$ is representable and hence, $\mathfrak{A}$ is representable.
\end{proof}
\end{Theorem}

As in the theory of cylindric algebras [2], a dimension-complemented $SA_\alpha$ is defined as follows:

\begin{Definition}
An $SA_\alpha$ $\mathfrak{A}$ is dimension-complemented if for all finite $X\subseteq A$, 
\[\alpha-\bigcup\{\Delta x:x \in X\}\]
 is infinite. 
\end{Definition}

\begin{Theorem}
Every dimension-complemented $SA_\alpha$ is embeddable in a dimension-complemented  $SA_\alpha$ in which elements are distinguished.
\begin{proof}

The proof is obtained from the proof of Theorem 5.3 and 5.4 of [1] if we replace ``locally finite'' wherever it occurs by ``dimension-complemented.''

\end{proof}
\end{Theorem}

In the theory of cylindrical algebras, every dimension-complemented cylindrical algebra is representable [3]. The following is a theorem analogous to this.

\begin{Theorem}
Every dimension-complemented $SA_\alpha$ is representable.
\begin{proof}
This follow immediately from Theorems 4.3 and 4.1.

\end{proof}
\end{Theorem}

\section{Neat embeddings}

\begin{Definition}
Let $\mathfrak{A}=\langle A, *_\kappa,v_\kappa\rangle_{\kappa <\alpha}$ be an $SA_\alpha$ and $\beta \leq \alpha$. The $\beta$-reduct of $\mathfrak{A}$ is the algebra $\mathfrak{A}/\beta=\langle A, *_\kappa,v_\kappa\rangle_{\kappa <\beta}$.
\end{Definition}

Clearly, $\mathfrak{A}/\beta$ is an $SA_\beta$ if $\mathfrak{A}$ is an $SA_\alpha$.

\begin{Theorem}
If $\mathfrak{A}$ is a $RSA_\alpha$ and $\beta \leq \alpha$, then $\mathfrak{A}/\beta$ is an $RSA_\beta$.
\begin{proof}
Let $\mathfrak{A}$ be an  $RSA_\alpha$. Then $\mathfrak{A}$ is isomorphic to a subalgebra of an $FSA_\alpha$ $\mathfrak{A'}$. $\mathfrak{A'}$ is a subalgebra of $\mathfrak{B}$, the full $FSA_\alpha$ with base $U$. Therefore, $\mathfrak{A'}/\beta$ is embeddable in $\mathfrak{B}/\beta$. Theorem 3.2 can be easily modified to show that elements are distinguished in $\mathfrak{B}/\beta$ and hence, by Theorem 4.1, $\mathfrak{A'}/\beta$ is an $RSA_\alpha$. $\mathfrak{A}/\beta$ is isomorphic to a subalgebra of $\mathfrak{A'}/\beta$. Therefore, $\mathfrak{A}/\beta$ is an $RSA_\beta$.

\end{proof}
\end{Theorem}

\begin{Definition}
Let $\mathfrak{A}$ be an  $SA_\alpha$ and $\beta \leq \alpha$. Let $\mathfrak{A}^{(\beta)}=\langle A',*_\kappa,v_\kappa \rangle_{\kappa < \beta}$ where $A'=\{s:s \in A \textrm{ and } \Delta x \subseteq \beta\}$. An  $SA_\alpha$ $\mathfrak{B}$ is $\beta$-neatly embeddable in $\mathfrak{A}$ if $\mathfrak{B}$ is isomorphic to a subalgebra of $\mathfrak{A}^{(\beta)}$.
\end{Definition}

As in the theory of cylindric algebras, necessary and sufficient conditions for representability of $SA_\alpha$'s can be given in terms of neat embeddings [3] .

\begin{Theorem}
Let $\mathfrak{A}$ be an $SA_\alpha$. The following statements are equivalent:

(i) $\mathfrak{A}$ is an $RSA_\alpha$.

(ii) For all $\kappa < \omega$, $\mathfrak{A}$ is $\alpha$-neatly embeddable in some $SA_{\alpha+\kappa}$.

(iii) $\mathfrak{A}$ is $\alpha$-neatly embeddable in some $SA_{\alpha+\omega}$.

\begin{proof}
To show that (i) implies (ii), it suffices to show that if $\mathfrak{A}$ is an $FSA_\alpha$ and $\kappa < \omega$, then $\mathfrak{A}$ is $\alpha$-neatly embeddable in some $FSA_{\alpha+\kappa}$. Let $\mathfrak{A}$ be an $FSA_\alpha$ with base $U$ and let $\kappa < \omega$. Define a function $\phi:A \to U^{U^{\alpha+\kappa}}$ by $\phi(f)(s)=f(s/\alpha)$ where $f \in A$, $s \in U^{\alpha + \kappa}$ and $s/\alpha$ is the restriction of $s$ to $\alpha$. It is easy to verify that $\phi$ maps $\mathfrak{A}$ isomorphically onto a subalgebra of $\mathfrak{B}^{(\alpha)}$ where $\mathfrak{B}$ is the full $FSA_{\alpha+\kappa}$ with base $U$.

Assume (ii). For $\kappa < \omega$, let $\phi_\kappa$ map $\mathfrak{A}$ isomorphically onto a subalgebra of $\mathfrak{B}_\kappa^{(\alpha)}$ where $\mathfrak{B}_\kappa=\langle B_\kappa,*_\lambda^\kappa,v_\lambda^\kappa\rangle_{\lambda<\alpha+\kappa}$ is an $SA_{\alpha+\kappa}$. Define algebras $\mathfrak{C_\kappa}= \langle B_\kappa,*_\lambda^\kappa,v_\lambda^\kappa\rangle_{\lambda<\alpha+\omega}$ where $x*_\lambda ^\kappa y=y$ for $x, y \in B_\kappa$ and $\alpha+\kappa \leq \lambda<\alpha+\omega$ and $v_\lambda^\kappa$ is any member of $B_\kappa$ for $\alpha+\kappa \leq \lambda<\alpha+\omega$. Let $\mathcal{F}$ be a non-principal ultafilter on $\omega$ and let $\mathfrak{C}=\prod \langle \mathfrak{C}_\kappa:\kappa<\omega\rangle/\mathcal{F}$.

We show that $\mathfrak{C}$ is an $SA_{\alpha+\omega}$. If an instance of an axiom for $SA_{\alpha+\omega}$'s has maximum subscript $\gamma$, then for all $\kappa$ with $\gamma<\alpha+\kappa$,this axiom holds in $\mathfrak{B}_\kappa$ and hence in $\mathfrak{C}_\kappa$ where $\gamma<\alpha+\kappa$. Hence, it holds in all but finitely many of the $\mathfrak{C}_\kappa$ and therefore, in  $\mathfrak{C}$.

Define a function $\phi: A \to C$ by $\phi(a)=[\langle \phi_\kappa(a):\kappa<\omega\rangle]$ where $a \in A$. It is easily verified that $\phi$ is an isomorphism from $\mathfrak{A}$ into a subalgebra of $\mathfrak{C}^{(\alpha)}$. Hence, we have (iii).

To show that (iii) implies (i), assume that $\mathfrak{A}$ is mapped isomorphically by $\phi$ onto a subalgebra of $\mathfrak{B}^{(\alpha)}$ where $\mathfrak{B}$ is a $SA_{\alpha+\omega}$. Using $\Delta(x*_\kappa y)\subseteq \Delta x \cup (\Delta y-\{\kappa\})$ [1], it can be shown that the subalgebra $\mathfrak{A}'$ generated by $ \{\phi(a):a \in A\}$ in $\mathfrak{B}$ is a dimension-complemented $SA_{\alpha+\omega}$. By Theorem 4.4, $\mathfrak{A'}$ is representable. $\mathfrak{A}$ is isomorphic to a subalgebra of $\mathfrak{A'}/\alpha$. By Theorem 5.2, $\mathfrak{A'}/\alpha$ is representable and hence, $\mathfrak{A}$ is representable.
\end{proof}
\end{Theorem}

\begin{center}
References
\end{center}

[1] N. \textsc{Feldman}, \textit{Representation of polynomial substitution algebras}, \textit{\textbf{Journal of symbolic logic}}, vol. 47 (1982), pp. 481-492.

[2] L.\textsc{Henkin}, J.D.  \textsc{Monk}, and A. \textsc{Tarski},  \textit{\textbf{Cylindric algebras, part I}}, North-Holland Publishing Company, Amsterdam,1971.

[3] L.\textsc{Henkin}, J.D.  \textsc{Monk}, and A. \textsc{Tarski},  \textit{\textbf{Cylindric algebras, part II}}, North-Holland Publishing Company, Amsterdam,1985.

[4] J. \L \textsc{o\'s},  \textit{Quelques remarques, th\'eor\`emes et probl\`emes sur les
              classes d\'efinissables d'alg\`ebres},   \textit{\textbf{Mathematical interpretations of formal systems}}, North Holland Publishing Company., Ansterdam, 1955 pp. 98-113.

\end{document}